\input amstex
\documentstyle {amsppt}

\advance\voffset  by -1.0cm
\NoBlackBoxes
\magnification=\magstep1
\hsize=17truecm
\vsize=24.2truecm
\voffset=0.5truecm
\document

\define \De {\Delta}
\define \cl {\Cal L_{S^1}}
\define \Ga {\Gamma}
\define \la {\lambda}
\define \pa {\partial}
\define \bR {\Bbb R}
\define \bC {\Bbb C}
\define \bP {\Bbb P}
\define \ga {\gamma}
\define \gas {\gamma: S^1\to \Bbb R^{2n}}
\define \pn {\Bbb RP^n}
\define \rtn {\Bbb R^{2n}}
\define \cg {\text{Cen}_\gamma}
\define \pnn {\Bbb RP^{2n}}
\define \part {\partial}

\topmatter
\title
On Young hulls of  convex curve in $\Bbb R^{2n}$\
\endtitle

\author V.~Sedykh and B.~Shapiro
\endauthor
\affil
Department of Mathematics, Moscow State University of Technology "Stankin",
Moscow 101472, Russia,
{\tt sedykh\@ium.ips.ras.ru},\\
 Department of Mathematics, University of Stockholm, S-10691, Sweden, {\tt
shapiro\@matematik.su.se}\
\endaffil
\rightheadtext {On Young hulls of convex curve}

\abstract
For a convex curve in an even-dimensional affine space we introduce a
series of convex domains (called Young hulls)   describe their structure
and
give a formulas fo the volume of the biggest of these domains.
\endabstract

\endsubjclass
\endkeywords

\endtopmatter

\heading
Introduction
\endheading

A simple smooth curve in $\Bbb R^m$ is called {\it convex} if the total
multiplicity of its
intersection with any affine hyperplane does not exceed $m$. Note that each
convex curve is
{\it nondegenerate}, i.e. has a nondegenerate osculating $m$-frame at each
point.

Nondegenerate curves in projective spaces have naturally arisen in many
classical geometrical
problems  as well as in problems related to linear ordinary differential
equations,
see e.g. \cite {KN,Shr,Li}. The global topological properties of the spaces of
such curves were used
 in the enumeration of symplectic leaves of the  Gelfand-Dickey bracket, \cite
{KSh}.

The set of all nondegenerate curves (say, with fixed osculating flags at
both endpoints) contains the
subset of convex curves which minimize the maximal
possible number of intersection
points with hyperplanes.
Convex curves correspond to a special
class of linear
ODE called disconjugate. Their properties
were studied, for example, in the classical papers \cite {Ga,Kr,Po} in
 connection with  Sturmian
theory and functional analysis. Recently V.~Arnold proposed a
generalization of Sturmian theory
 which  in the simplest treated case deals with the estimation of the
number of flattening points of
a curve in $\Bbb R^{2n+1}$ which projects on a curve close to the
generalized ellipse (1),
see \cite {Ar}. (For the reader interested in the general theory of
disconjugate linear ODE and convex curves we recommend the survey books
\cite {KN} and \cite {Co}.) In
the present paper we associate to each
closed convex curve in $\bR^{2n}$ a number of convex domains
and study some of their properties.

Note that a convex curve  in an odd-dimensional space can not be  closed. A
simple example of a closed convex curve
in $\Bbb R^{2n}$ is the {\it normalized generalized ellipse}
$$
(\sin t,\cos t,1/2\sin {2t},1/2\cos {2t},\dots ,1/n\sin {nt}, 1/n\cos
{nt}).\tag1
$$

{\smc Definition.} A tangent hyperplane to a curve is called a {\it support
hyperplane} if
the curve lies on one side w.r.t. this hyperplane.

Let $\ga : S^1\to \rtn$ be a simple smooth closed convex curve in $\Bbb
R^{2n}$. Then for any
point $t$ of $\ga$ there exists a support hyperplane passing through this point. For instance,
the {\it osculating} hyperplane to $\ga$ at $t$ is a support hyperplane
(the multiplicity of the
intersection of this hyperplane with $\ga$ at $t$ is equal to $2n$).

{\smc Definition.} The {\it convex hull} of a curve in an affine space is
the intersection of
all closed half-spaces contained this curve.

I.~Schoenberg [Sch] proved that the convex hull $\Cal {CH}_{\ga}$ of a
convex curve $\ga$ in $\rtn$ is the
intersection of a family of closed half-spaces determined by support
 hyperplanes
which are tangent to $\ga$ at $n$ distinct points. Moreover,  the Euclidean
volume $Vol(\Cal {CH}_{\ga})$
of its convex hull can be expressed by the integral
$$
Vol(\Cal {CH}_\ga)=\pm \frac {1}{n!(2n)!}\int _{\Cal T^n}\det
\bigl[\ga (t_1),\dots ,\ga (t_n),\ga \prime(t_1),\dots ,\ga
\prime(t_n)\bigr]dt_1\dots dt_n,\tag2
$$
where $ \Cal T^n=(S^1)^n$ denotes the $n$-dimensional torus.

Note that the support hyperplanes considered by Schoenberg intersect  $\ga$
with
the maximal possible multiplicity $2n$. But there exist other types of
support hyperplanes
intersecting $\ga$ with such a total multiplicity. These different types
are enumerated by the
Young diagrams of area $n$.

{\smc Definition.} A support hyperplane to a curve $\ga$ is said to be a
{\it hyperplane of
a given Young type} $\mu=(k_1,...,k_r)$ or just {\it a $\mu$-hyperplane} if
it is tangent to $\ga$ at $r$
pairwise different points
$\ga(t_1),...,\ga(t_r)$ with multiplicities $2k_1,...,2k_r$ resp.

Any Young diagram of  area at most $n$ defines in the above mentioned way the
family  of support hyperplanes to a convex curve $\ga$. Conversely, any
support hyperplane to a curve $\ga$
is a support hyperplane with some Young diagram of  area at most $n$.

{\smc Remark.} Note that for a convex curve $\ga: S^1\to \bR^{2n}$ and an
$r$-tuple of points $(t_1,...,t_r)$ the hyperplane tangent
to $\ga$ at $\ga(t_1),...,\ga(t_r)$ with the multiplicities
 $2k_1,...,2k_r$ resp. where $\sum k_i=n$ is automatically
a support hyperplane. This explains many simplifications which take place
  for convex curves in comparison with the general case.

{\smc Definition.} The {\it Young $(k_1,...,k_r)$-hull} of a convex curve
in $\ga : S^1\to \bR^{2n}$ is
the intersection of all closed half-spaces containing $\ga$  and defined by
the support
hyperplanes of a given Young type $(k_1,...,k_r)$.

Denote by $\Cal {YH}_{\ga}(k_1,...,k_r)$ the Young $(k_1,...,k_r)$-hull of
a curve $\ga$.

{\smc Remark.} If the
area of the Young diagram $(k_1,...,k_r)$ is less than $n$, then $\Cal
{YH}_{\ga}(k_1,...,k_r)=\Cal {YH}_{\ga}
(1,...,1,k_1,...,k_r)$ where the area of Young diagram
$(1,...,1,k_1,...,k_r)$ is equal to $n$.

Consider now all Young hulls of a curve $\ga$ defined by all Young diagrams
of  area $n$.
These hulls are enclosed according to the lexicographic
 order on Young diagrams. Namely,
$\Cal {YH}_{\ga}(k_1,...,k_r) \subset \Cal
{YH}_{\ga}(k_1^\prime,...,k_{r^\prime}^\prime)$ iff
$(k_1,...,k_r)$ is smaller than $(k_1^\prime,...,k_{r^\prime}^\prime)$ in
the standard
lexicographic order.

Thus, for any Young diagram $(k_1,...,k_r)$ of  area $n$, we have
$$
\Cal {YH}_{\ga}(1^n) \subset \Cal {YH}_{\ga}(k_1,\dots ,k_r) \subset \Cal
{YH}_{\ga}(n),
$$
where $(1^n)\equiv (1,...,1)$ ($n$ times). The smalest Young hull $\Cal
{YH}_{\ga}(1^n)$ of a
curve $\ga$ is its convex hull $\Cal {CH}_{\ga}$. The biggest Young hull
$\Cal {YH}_{\ga}(n)$ also denote by
$\Cal {EH}_{\ga}$ of a curve $\ga$ is called its {\it elliptic hull}.

The main result of this paper is a characterization of the  Young hulls of
convex curves
in even-dimensional spaces  and their duals as the convex hulls of some
special'varieties' as well as a formula
for the volume of the elliptic hull $\Cal {EH}_{\ga}$ similar to (2). As an
example of the obtained results, let us describe the structure of the
elliptic hull $\Cal {EH}_\ga$.

Let $L(t)$ be a $(2n-2)$-dimensional affine subspace in $\Bbb R^{2n}$ passing
through the point $\ga(t)$  and spanned by the vectors $\ga^\prime(t),\dots ,
\ga ^{2(n-1 )}(t)$. Then we have a map $\Ga_\ga: (\Cal T^n\setminus
Diag)\to\rtn$ sending an
$n$-tuple of pairwise different points $(t_1,\dots , t_n)\in \Cal T^n$ to
the intersection
$L(t_1)\cap ...\cap L(t_{n-1})\cap L(t_n)$. Obviously, $\Ga_\ga$
is invariant w.r.t. the action of the symmetric group
$\frak S_n$ on $\Cal T^n$ by permutation of coordinates.

The map $\Ga$ extends continuously to the
$\frak S_n$-invariant map $\Ga_\ga: \Cal T^n\to\rtn$ (in particular,
$\Ga_\ga(t,t,...,t)$ coincides with $\ga(t)$).

{\smc Proposition.}
The elliptic hull $\Cal {EH}_\ga$ of a curve
$\ga$ is the convex hull of the image $\Ga_\ga (\Cal T^n)=$ $\Ga_\ga
 (\Cal T^n/\frak S_n)$ and its Euclidean volume can be
expressed by the integral
$$\multline
Vol(\Cal {EH}_\ga)=\pm \frac {1}{(2n)!}\int _{\Cal T^n}\det
\bigl [\Ga_\ga(t_1,t_1,...,t_1),\Ga_\ga(t_1,t_2,t_2,...,t_2),\dots
,\Ga_\ga(t_1,...,t_n),\\
\frac{\part \Ga_\ga(t_1,...,t_n)}{\part t_1},\dots ,
\frac{\part \Ga_\ga(t_1,...,t_n)}{\part t_n}\bigr ] dt_1...dt_n.
\endmultline\tag3
$$

Recall that using  (2), Schoenberg has  proved
the following interesting isoperimetric
inequality for convex curves
$$
l^{2n}_\ga \geq (2\pi n)^nn!(2n)!Vol(\Cal {CH}_\ga)\tag4
$$
where $l_\ga$ denotes the Eucledian length of $\ga$. Moreover, he has shown
that this estimation is
sharp and the equality holds only for the curve (1) up to a parallel shift,
a homothety
and an orthogonal transformation of $\bR^{2n}$.

The main motivation of the present study was an attempt to generalize
Schoenberg's inequality for the
case of elliptic (and more general Young) hulls. Unfortunately,  at the
present moment
this goal still remains unachieved even for the elliptic hulls. We hope to
return to this problem
in the future.

{\smc Conjecture.} An estimate of the form
$$l_\ga^{2n}\ge \kappa_\mu Vol(\Cal{YH}_\ga(\mu))$$
( where $\kappa_\mu$ is some constant depending on $\mu$
 only) is valid for all Young hulls of a curve $\ga$.
Moreover, the equality holds only for the curve (1) (considered up to
a homothety and the group of rigid motions of $\bR^{2n}$).

The structure of the paper is as follows.  In \S 1 we prove
some basic facts about $\Cal {YH}_\ga(\mu)$, in particular,
find a special ruled hypersurface bounding $\Cal {YH}_\ga(\mu)$ and
describe the minimal convex subset spanning $\Cal {YH}_\ga(\mu)$ and its
dual as its convex hull.  In \S 2
we present an integral formula for $Vol(\Cal {EH}_\ga)$.
Finally, in \S 3 we calculate explicitly the volume of the elliptic hull of
the generalized normalized ellipse in $\bR^4$.

 {\smc Acknowledgements.} The
first author is grateful to the Swedish Natural
Science Research Council for the financial support  during his stay at the
Department of Mathematics of the University of Stockholm where this project
was started. The second author wants to acknowledge the hospitality and
nice research atmosphere of the Max-Planck Institute during the
preparation of the manuscript.

\heading \S 1. Generalities on  convex curves \endheading

\subheading {Boundness}

We start with the following simple proposition.

{\smc 1.1. Lemma.} For a convex $\ga: S^1\to \rtn$ every Young hull
$\Cal {YH}_\ga(\mu)$ is a  bounded convex domain.

{\smc Proof.} It suffices to show that $\Cal{EH}_\ga$ is bounded since it
contains all other Young hulls. One has to show that for every point $p$
lying outside a sufficiently
big ball containing $\ga$ there exists an osculating  hyperplane to $\ga$
passing through $p$. This is equivalent
to finding a hyperplane $H_p$ passing through $p$ and tangent to $\ga$.
Taking a pont $p$ as above we can find a hyperplane  $\tilde H_p$ through
$p$ not intersecting $\ga$ at all. Now
fixing some codimension 2 subspace $\tilde h_p\subset \tilde H_p$ through
$p$ we can take the pencil of hyperplanes
 containing $\tilde h_p$. Since $\tilde H_p$ does not intersect $\ga$ this
family of hyperplanes contains a
hyperplane tangent to $\ga$.

\qed

\subheading {$\mu$-discriminants}

{\smc Definition.} Given $\mu=(k_1,...,k_r),\,\sum k_i=m$ we call a  curve $\ga: S^1\to \Bbb R^m$  {\it $\mu$-generic}
if  for any $r$-ruple of pairwise different points $(t_1,...,t_r)$
 {\it the spanning affine space $L_\mu(t_1,...,t_r)$} passing through
$\ga(t_1),\ga(t_2),...,
\ga(t_r)$ and containing
$\ga^\prime(t_1),...,\ga^{(k_1-2)}(t_1),\ga^\prime (t_2),...,$
$\ga^{(k_2-2)}(t_2),...,$ $\ga^{(k_r-2)}(t_r)$ has codimension $r+1$.

One can easily check that $\mu$-generic curves form an open dense subset in the
space of all maps $S^1\to \bR^m$.

{\smc Definition.} Given a Young diagram $\mu=(k_1\ge k_2\ge ...\ge
k_r),\;\sum_ik_i=m$
we call by {\it the $\mu$-discriminant $D_\ga(\mu)$} of a $\mu$-generic
curve $\ga:S^1\to \bR^m$ the
ruled hypersurface obtained as the closure of the union of all spanning
affine subspaces  $L_\mu(t_1,...,t_r)$.

{\smc Remark.} If $\mu=(m)$ then we get the standard discriminant $D_\ga$
(also called
the swallowtail or the front of $\ga$), i.e. the hypersurface ruled out by
the family of $(m-2)$-dimensional osculating subspaces.
 In a neighborhood of a generic point of $\ga$ (where all
$m$ derivatives
$\ga^\prime(t),\dots ,\ga ^{(m)}(t)$ are linearly independent) a germ of $D$
is diffeomorphic to
a germ of the usual discriminant consisting of all polynomials of degree
$m$ with multiple zeros.
If $\ga$ is globally nondegenerate, i.e. $\ga^\prime(t),\dots ,\ga
^{(m)}(t)$ are linearly
independent for any $t\in S^1$ then for any point $p\in
\Bbb R^m\setminus
D_\ga$ and any hyperplane $\tilde H_p$ through $p$
the local multiplicity of intersection of $\tilde H_p$ and $\ga$ does not
exceed $2$.

{\smc 1.2. Proposition.\;} Given a convex curve $\ga:S^1\to \bR^{2n}$ and a
diagram $\mu=(k_1,...,k_r),$
$\sum k_i=n$ one gets that $\Cal {YH}_\ga(\mu)$ coincides with a convex
connected component of the
complement $\bR^{2n}\setminus D_\ga({2\mu})$ where $2\mu=(2k_1,...,2k_r)$
is the 'doubled' diagram of area $2n$.

{\smc Proof.} Proposition follows immediately from lemmas 1.3-1.4.

\qed

{\smc Definition.} The intersection of the Young $\mu$-hull $\Cal
{YH}_\ga(\mu)$  of a curve $\ga$ with a given support hyperplane is called
the {\it characteristic set} of this hyperplane.

The characteristic set of a support hyperplane is convex (as the
intersection of 2 convex sets).

{\smc 1.3. Lemma.} The boundary of the Young $\mu$-hull  $\Cal {YH}_\ga(\mu)$
of $\ga$ is the closure of the union of characteristic sets of all support
$\mu$-hyperplanes to $\ga$.

{\smc Proof.} The characteristic set $X$ of a given hyperplane $\pi$ can not
contain interior points of $Y=\Cal {YH}_\ga(\mu)$ since $Y$ lies on one
side of $\pi$. Hence, $X$ belongs to the boundary $\part Y$.

Conversely, let $p\in \part Y$. Then there exists a sequence of points
$p_i\notin Y$ and $\mu$-hyperplanes $\pi_i,\, i=1,...,$ such that $p_i$
and $Y$ lie on the different sides of $\pi_i$ and $\lim_{i\to \infty}p_i=p$.

Taking subsequences of hyperplanes we can assume that there exists
$\lim_{i\to \infty}\pi_i=\pi$. Then the hyperplane $\pi$ is a support
hyperplane to $\ga$. The point $p$ and the set $Y$ have to lie on different
sides of
$\pi$. Hence, $p\in \pi\cap Y$ and $p\in Y$.

\qed

{\smc 1.4. Lemma.} The characteristic set of a support $\mu$-hyperplane
tangent to
$\ga$ at $\ga(t_1),...,\ga(t_r)$ belongs to the spanning affine subspace
$L_\mu(t_1,...,t_r)$.

{\smc Proof.} Let $D[u_1,...,u_r,v]$ be the determinant of order $2n$
of the matrix consisting of
$\ga(t_2)-\ga(t_1),...,\ga(t_r)-\ga(t_1),\ga^\prime(t_1),...,
\ga^{(2k_1-2)}(t_1),...,\ga(t_r),...,\ga^{(2k_r-2)}(t_r),u_1,...,u_r,v$.
Then for any small $\De_1,...,\De_r$ the support $\mu$-hyperplane tangent
to $\ga$ at $\ga(t_1+\De_1),...,\ga(t_r+\De_r)$ defined by the equation
$F(x,\De_1,...,\De_r)=0$ where $x\in \rtn$ and

$$\align &F(x,\De_1,...,\De_r)=
D[\ga^{(2k_1-1)}(t_1),...,\ga^{(2k_r-1)}(t_r),x-\ga(t_1)]\\
&+D[\ga^{(2k_1)}(t_1),\ga^{(2k_2-1)}(t_2),...,\ga^{(2k_r-1)}(t_r),x-\ga(t_1)]
\De_1+...\\
&+D[\ga^{(2k_1-1)}(t_1),...,\ga^{(2k_{r-1}-1)}(t_{r-1}),
\ga^{(2k_r-1)}(t_r),x-\ga(t_1)]\De_r+...\endalign$$

If $x$ lies on the support hyperplane $\pi$ tangent to $\ga$ at
$\ga(t_1),...,\ga(t_r)$ and does not belong to the subspace
$L_\mu(t_1,...,t_r)$ then $F(x,0,...,0)\neq 0$ and there exists
$i\in [1,...,r]$ such that $\part F/\part\De_i(x,0,...,0)\neq 0$
(since $\ga$ is convex and the area of $\mu$ equals $n$).
Hence one can find $\De_i^-<0$ and $\De_i^+>0$ such that
$$F(x,0,...,0,\De_i^-,0,...,0)F(x,0,...,0,\De^+_i,0,...,0)<0.$$
Therefore $x$ can not belong to the characteristic set of $\pi$.

\qed

{\smc Remark.} Appearently for a  convex $\ga:S^1\to \rtn$ the complement
$\rtn\setminus D_{\ga}(2\mu)$ contains a unique convex connected component.

{\smc Example 1} (useful although the source is $\Bbb R$ instead of $S^1$).
Consider the usual
rational normal curve $\xi_m: \Bbb R\to \Bbb R^m,\; t\to (mt,\binom m2
t^2,\dots , \binom m2 t^{m-1},
t^m)$. Identifying $\Bbb R^m$ with the space $\frak {Pol}_m$ of all monic
polynomials $x^m+a_1x^{m-1}+\dots
+a_m$ in the obvious way one identifies the curve $\xi_m$ with the
$1$-parameter family of
polynomials $(x+t)^m$.
 The elliptic hull $\Cal {EH}_{\xi_m}$ consists of all strictly elliptic
polynomials, i.e. polynomials with no real
zeros. ($\Cal{EH}_{\xi_m}$ is nonempty only for even $m$.) The last observation
motivates the term 'elliptic hull' of $\ga$.

\subheading {Dual domain and convex skeleton}

For any nondegenerate $\ga:S^1\to \Bbb R^m$ one can define its dual
$\ga^*:S^1\to (\bP^{m})^*$ considering
$\Bbb R^m$ as the standard affine chart in $\bP^m$.

{\smc Definition.\;} The curve $\ga^*:S^1\to (\bP^{m})^*$ such that
$\ga^*(t)$ is the hyperplane
spanned by $\ga^\prime(t),\dots , \ga^{(m-1)}(t)$ is called {\it the dual
curve} to the nondegenerate
curve $\ga$.

{\smc Remark.\;} The dual to a convex curve is a convex curve, see e.g.
\cite {Ar}.

Given a convex closed domain $\Omega\subset\Bbb R^m\subset \bP^m$ (here
'convexity'
means that with any
two points $p_1,p_2$ the set $\Omega$ contains the whole segment
$(p_1,p_2)\subset \bR^m$)
 one defines the dual convex domain as follows.

{\smc Definition.} {\it The dual convex domain} $\Omega^*\subset (\bP^m)^*$
is the
closure of the set of all hyperplanes not intersecting $\Omega$.

{\smc Remark.} One can easily check that $\Omega^*$ lies in
some affine chart of $(\bP^m)^*$ and is convex there in the usual sense.
Namely,
fixing a point $p$ in the interior of $\Omega$ we get the corresponding
   hyperplane $H_p^*\subset
(\bP^m)^*$. Obviously, $\Omega^*\subset (\bP^m)^*\setminus H_p^*$.

Now we will describe the  minimal subset the convex hull of which coincides
with  the domain dual to some Young hull $\Cal {YH}_\ga(\mu)$.

{\smc Definition.} Given a convex $\ga: S^1\to \rtn$ and
 an $n$-tuple of pairwise different points $(t_1,t_2,\dots , t_n)\in \Cal
T^n\setminus Diag$ we call by {\it the associated map} $\Ga_\ga: \Cal T^n
\setminus Diag\to \rtn$ the map sending $(t_1,t_2,\dots ,t_n)$ to the
intersection point $\bigcap_{i=1}^nL(t_i)$ where
$L(t)$  denotes the $(2n-2)$-dimensional subspace passing through $\ga(t)$
and spanned by $\ga^\prime(t),\dots ,$ $ \ga ^{(2n-2)}(t)$.

Convexity of $\ga$ provides that this intersection is a point for any
choice of $(t_1,t_2,\dots t_n)\in \Cal T^n\setminus Diag$.
(Note that $L(t)$ are subspaces ruling out the swallowtail $D_\ga$.)

{\smc Definition.\;} {\it By the convex skeleton $\bold {CSk}_\ga$} of a convex
$\gas$ we denote the
closure of the image of the associated map   $\overline{\bigcup_{t_1,\dots
,t_n\in \Cal T^n\setminus Diag}\Ga_\ga(t_1,t_2,\dots
,t_n)}$.

{\smc 1.5. Lemma.} The associated map
$\Ga_\ga: (\Cal T^n\setminus Diag)\to\rtn$    extends
continuously to the $\frak S_n$-invariant map $\Ga_\ga: \Cal T^n\to\rtn$
where as before $\frak S_n$ is the symmetric group
acting by permutation of coordinates of $\Cal T^n$.

{\smc Proof.} One can easily see that convexity of $\ga: S^1\to \rtn$
implies (and is in
fact equivalent) to the following property. For any $r$-tuple of
points $\ga(t_1),...,\ga(t_r)$ there exists and unique hyperplane passing
through these points and tangent to $\ga$ with the total multiplicity $2n$
distributed in any way between  $\ga(t_1),...,\ga(t_r)$.
In our case using the dual space one gets that the intersection of
subspaces $\bigcap_{i=1}^nL(t_i)$ corresponds to the hyperplane
$H_{t_1,...,t_n}\subset (\bP^{2n})^*$ tangent to the dual curve $\ga^*$ at
$\ga^*(t_1),...,\ga^*(t_n)$ with the multiplicity 2 at each point. Note
that if $\ga$ is convex
then $\ga^*$ is also convex. The closure of the set
$\bigcup_{(t_1,...,t_n)\in \Cal T^n\setminus Diag} H(t_1,...,t_n)$,
obviously, is the closed set of hyperplanes tangent of $\gamma^*$ with the
even local multiplicity at every point. Going back to $\rtn$ we get the
required result.

\qed

{\smc Example 2.} The convex skeleton $\bold {CSk}_{\xi_{2n}}$ in the
example 1 is the set of
all polynomials of the
form $(x-t_1)^2...(x-t_n)^2$.

{\smc 1.6.
Corollary.\;} The map $\Ga_\ga: \Cal T^n/\frak S_n\to \rtn$ sends $\Cal
T^n/\frak S_n$
homeomorphically onto $\bold {CSk}_\ga$.

{\smc Proof.} It suffices to show that $\Ga_\ga:\Cal T^n/\frak S_n\to\rtn$ is
a 1-1-map. This again follows from the fact that for a convex $\ga^*$
and a pair $(t_1^1,...,t_{r^1}^1,k_1^1,...,
k_{r^1}^1)$ and $(t_1^{2},...,t^{2}_{r^{2}},
k_1^{2},...,k_{r^{2}}^{2})$,
$\sum k_i^1=n$ and $\sum k_i^{2}=n$ the hyperplanes $H_1$ tangent to $\ga$ at
$(t_1^1,...,t_{r^1}^1)$ with the
multiplicities $(2k_1^1,...,2k_{r^1}^1)$ and $H_2$ tangent to $\ga$ at
$(t_1^{2},...,t^{2}_{r^2})$ with the multiplicities
$(2k_1^{\prime\prime},...,2k^{\prime\prime}_{r^{\prime\prime}})$   are
different unless $(t_1^1,...,t_{r^1}^1,k_1^1,...,
k_{r^1}^1)=(t_1^{2},...,t^{2}_{r^{2}},
k_1^{2},...,k^{2}_{r^{2}})$.

\qed

The above homeomorphism provides the following natural stratification of
$\bold {CSk}_\ga$ into strata enumerated by the Young diagrams of area $n$.
Fixing a partition $\mu=(k_1\geq k_2\geq \dots \geq k_r),\;\sum k_i=n$
let us denote by $\bold {CSk}_\ga (\mu)$  the closure of the set of all
intersections
$\bigcap_{i=1}^rL_{k_i}(t_i)$ where $(t_1,...,t_r)$ are pairwise different
and $L_{k_i}(t_i)$
is the subspace passing through $\ga(t)$ and spanned by
$\ga^\prime(t_i),...,\ga^{2(n-k_i)}(t_i)$.
Obviously,
$\dim \bold {CSk}_\ga (\mu)=r$ and different strata $\bold {CSk}_\ga (\mu)$
form the stratification
of $\bold {CSk}_\ga$ with the adjacency  coinciding with the reverse
lexicographic
partial order on the set of all Young diagrams of area $n$.
Namely, a stratum
$\bold {CSk}_\ga(k_1,...,k_{r_1})$
is adjacent to $\bold {CSk}_\ga(k_1^\prime,...,k_{r_2}^\prime)$  iff the
partition
$(k_1,...,k_{r_1})$ is smaller  than $(k_1^\prime,...
,k_{r_2}^\prime)$ in the  lexicographic order.

{\smc 1.7. Lemma.} The convex skeleton $\bold {CSk}_\ga$ (and each stratum
 $\bold {CSk}_\ga(k_1,..,k_r)$) is a weakly convex set, i.e. it
lies on the boundary of its convex hull.

{\smc Proof.} We will prove by induction on $n$ that $\bold {CSk}_\ga$ lies
on the boundary of $\Cal{EH}_\ga$ and therefore is weakly convex.
Base of induction, $n=1$. In this case
$\ga$ is a convex curve on $\bR^2$ and the statement is obvious.

Step of induction. Let us assume
that the statement holds for $n\le N$ and consider a convex curve $\ga$ in
$\bR^{2N+2}$. Let us show that all intersections $\bigcap_{i=1}^{N+1}
L(t_i)$ belong to
$\part \Cal {EH}_\ga$. Fixing  $t_1$ let us  take the following curve
$\ga^{t_1}\subset L(t_1)$. Namely, $\ga^{t_1}:S^1\to L(t_1)$ obtained as the
projection of $\ga$ along the family of  2-dimensional subspaces through
$\ga(t)$ containing
$\ga^\prime(t)$ and $\ga^{\prime\prime}(t)$ on $L(t_1)$.
One can show that $\ga^{t_1}$ is a smooth convex curve in $L(t_1)$. Moreover,
$\Cal{EH}_{\ga^{t_1}}$ coincides with $\Cal{EH}_\ga\cap L(t_1)$. This follows
from the fact that the intersection of an osculating subspace to $\ga$
restricted to $L(t_1)$ coincides with the osculating subspace to
$\ga^{t_1}$ at the corresponding point. Thus the closure of the union of
intersections
$\bigcap_{i=1}^N L(t_i)$ coincides with $\bold {CSk}_{\ga^{t_1}}$.
Since $\ga^{t_1}$ is convex in $\rtn$ everything is proved for $t_1$. Varying
$t_1$ we get the necessary statement.

\qed

Now we are in position to characterize the dual domains to $\Cal {YH}_\ga(\mu)$.

{\smc 1.8. Proposition.\;} For any convex $\ga:S^1\to \rtn$ and any diagram
$\mu$ of area $n$
the dual domain $\Cal {YH}_\ga^*(\mu)$ coincides with the convex hull of
the stratum $\bold {CSk}_{\ga^*}(\mu^*)$ where $\mu^*$ is the dual Young
diagram,
i.e. obtained by making rows into columns.

{\smc Proof.} By definition the boundary of the dual domain $\Cal
{YH}_\ga^*(\mu)$ consists of all support hyperplanes to $\Cal
{YH}_\ga(\mu)$.
One checks immediately that the union of all support $\mu$-hyperplanes to $\ga$
forms $\bold{CSk}_{\ga^*}(\mu^*)$. Fixing a point $p$ in the
interior of $\Cal {YH}_\ga(\mu)$ one gets that  $\Cal {YH}_\ga^*(\mu)\subset
(\bP^{2n})^*\setminus H_p^*$. For any other support hyperplane to $\ga$ its
halfspace
not containing $\ga$ lies in the union of halfspaces cut off by the
$\mu$-hyperplanes. This means that the hyperplane itself as the point in
the dual space lies  in the convex hull of
 $\bold {CSk}_{\ga^*}(\mu^*)$ in $(\bP^{2n})^*\setminus H_p^*$.

$\qed$

{\smc Observation.} The stratum $\bold {CSk}_{\ga}(\mu)$ is minimal, i.e.
none of its point $p$ lies in the convex hull of $\bold {CSk}_{\ga}(\mu)
\setminus p$.

{\smc Proof.} Interpreting a point $p\in \bold {CSk}_{\ga}(\mu)$ as a
support hyperplane of the  $\mu^*$-type to $\ga^*$ one can easily check that
in the case of convex $\ga$ a positive linear combination of 2 such hyperplanes
is again of $\mu^*$-type if and only if these hyperplanes coincide. (Using
interpretation with disconjugate linear ODE, see e.g. \cite {BSh} one can
 think of support hyperlanes of $\mu^*$-type as of nonnegative $2\pi$-periodic
solutions to a certain linear ODE having exactly $r$ zeros of multiplicity
 $2\mu_1^*,...,2\mu_r^*$ at some points $t_1,...,t_r$ on the period $[0,2\pi)$.
This interpretation makes the observation obvious.)

\subheading {$\Cal {YH}_\ga(\mu)$ as a convex hull}

Finally, we describe $\Cal {YH}_\ga(\mu)\;$ as the convex hull of some
special set in the initial space
$\rtn$. Let us return to the hypersurface $D_\ga(2\mu)$ ruled out by
codimension $r+1$ subspaces
$L_\mu(t_1,...,t_r)$.

{\smc Definition.}  We call the set of subspaces $L_\mu(t^1_1,...,t^1_r)$,
$L_\mu(t^2_1,...,t^2_r)$,
..., $L_\mu(t^p_1,...,t^p_r)$ {\it essential} if the intersection
$\bigcap_{i=1}^p L_\mu(t^i_1,...,t^i_r)$ 
is a point. This intersection point itself is called {\it essential} and,
finally, the closure of the union of
all essential points is called {\it the essential $\mu$-set
$\bold {Ess}_\ga(\mu)\;$} of $\ga$.

{\smc 1.9. Proposition.} For any convex $\ga:S^1\to \rtn$ and any diagram
$\mu$ of area $n$
the Young hull $\Cal {YH}_\ga(\mu)$ is the convex hull of $\bold
{Ess}_\ga(\mu)$.

In particular, one has

{\smc 1.10. Corollary.\;} For any convex $\gas$ the elliptic hull
${\Cal {EH}_\ga}$ is the convex hull of $\bold {CSk}_\ga$.

To prove statement 1.9. we need the following definitions.

{\smc Definition.} A support hyperplane $H$ to a convex domain $\Omega$ is
called  {\it flexible} if there exists a subspace $h\subset H$ of codimension 2
such that all hyperplanes containing $h$ and lying sufficiently close to $H$
are support hyperplanes. A support hyperplane which is not flexible is called
{\it rigid}.

{\smc Observation.} Obviously, the union of all rigid support hyperplanes
to $\Omega$ is the minimal weakly convex set spanning $\Omega^*$ as its
convex hull.

{\smc Proof of 1.10.} The boundary of $\Cal {YH}_\ga(\mu)$ consists of all
support hyperplanes to  $\Cal {YH}_\ga^*(\mu)$. By proposition 1.8. the
set  $\Cal {YH}_\ga^*(\mu)$ is the convex
hull of $\bold {CSk}_{\ga^*}(\mu^*)$. Therefore any hyperplane intersecting
$\bold {CSk}_{\ga^*}(\mu^*)$ at some points corresponds to the point in the
initial space contained in some intersection  $\bigcap_{i=1}^{\tilde
p}L_\mu(t_1^i,...,t_r^i)$. Assume now that $H_p$ is a support hyperplane to

 $\Cal {YH}_\ga^*(\mu)$ such that the corresponding intersection
 $\bigcap_{i=1}^{\tilde p}L_\mu(t_1^i,...,t_r^i)$ is at least 1 -dimensional.
Choosing  in this intersection any line through $p$ one gets a 1-parameter
family of supporting hyperplanes used in the definition of flexibility.
Thus $H_p$ is flexible and is not included in the minimal weakly convex set
of  $\Cal {YH}_\ga^*(\mu)$.  On the other side using induction similar to the
one from 1.7 one can show that each point in  $\bold {Ess}_\ga(\mu)$ lies in
 $\part\Cal {YH}_\ga(\mu).$

\qed

{\smc Example 3.} Any elliptic polynomial is a convex linear combination of
polynomials of the form $(x-t_1)^2...(x-t_n)^2$.

{\smc Example 4.} The intersection of $L _{(1^n)}(t_1,\dots ,t_n)$ with the
convex hull
$\Cal {CH}_{\ga}$ of a curve $\ga$ is the $(n-1)$-dimensional simplex with
vertices at points
$\ga (t_1),\dots ,\ga (t_n)$.

{\smc Example 5.} The intersection of $L_{(n)}(t_1)=L(t_1)$ with the
elliptic hull $\Cal {EH}_{\ga}$ of
a curve $\ga $ is the set bounded by the hypersurface in $L(t_1)$
consisting  of all intersection  points of
subspaces $L(t_1)\cap\dots \cap L(t_n)$ for all $n$-tuple of pairwise
different points
$(t_1,\dots , t_n)$ on $\ga$ with  $t_1$ fixed.

{\smc Example 6.} The intersection of $L_{(1,2)}(t_1,t_2)$ with the Young hull
$\Cal {YH}_{\ga}(1,2)$ of a curve $\ga$ in $\Bbb R^{6}$ is the set bounded
by the  surface consisting
of all intersection points with $L_{(1,2)}(t_3,t_4)$
where $(t_3,t_4)$ are varying.

\heading \S 2. Integral presentation of  volume of $\Cal {EH}_{\ga}$\endheading

$r+1\vert 2n$. In this case
essential set




Let as above $\Ga_\ga:\Cal T^n\to\rtn$ denote the $\frak S_n$-invariant
associated map parameterizing the convex skeleton $\bold {CSk}_\ga$. We use

the convention that $\Ga_\ga(t_1,...,t_j),\,j<n$ means that the last
variable $t_j$ is repeated $n-j+1$ times.

{\smc 2.1. Proposition.} The  family of simplices spanned by
$\Ga_\ga(t_1),\Ga_\ga (t_1,t_2),$ $...,$
$\Ga_\ga (t_1,...,t_n)$ where  $(t_1,...,t_n)\in \Cal T^n$ 'parametrizes'
the boundary  $\part \Cal {EH}_\ga$. Namely, each generic point
on $\part \Cal {EH}_\ga$ belongs to exactly $n!$ such simplices.

{\smc Proof.} Induction on $n$. For $n=1$ each point of a convex $\ga$
is passed exactly once when $t$ runs over $S^1$. Assume that the statement
 holds for $n\le N$ and consider a convex curve $\ga:S^1\to \bR^{2N+2}$.
A generic point $p\in \part\Cal {EH}_\ga$ belongs exactly to $N+1$
codimension 2 subspaces $L(t_1),...,L(t_{N+1})$. Take $L(t_1)$ and the
projected curve $\ga^{t_1}:S^1\to L(t_1)$ introduced in the proof of 1.7.
One has $p\in \part\Cal {EH}_{\ga^{t_1}}$ and the family of
$(N-1)$-diimensional
simplices for $\ga^{t_1}$ is obtained from the family
$\Ga_\ga(t_1),\Ga_\ga (t_1,t_2),$ $...,$
$\Ga_\ga (t_1,...,t_n)$ where  $(t_1,...,t_n)\in \Cal T^n$ by fixing $t_1$ and
varying $t_2,...,t_{N+1}$. Therefore by the inductive hypothesis in the
restricted family of simplices one covers $p$ exactly $N!$ times. Repeating
the argument for $t_2,...,t_{N+1}$ one gets that the whole family covers
$p$
exactly $(N+1)!$ times.

\qed

{\smc 2.2. Proposition.}
One has the following integral presentation
for $Vol(\Cal {EH}_\ga)$

$$\align Vol(\Cal {EH}_\ga)=&\pm \frac {1}{2n!}\int _{\Cal T^{n}}\det
[\Ga_\ga (t_1),\Ga_\ga (t_1,t_2),
...,\Ga_\ga (t_1,...,t_n),
\frac{\part \Ga_\ga (t_1,...,t_n)}{\part t_1},\\
& \frac{\part
\Ga_\ga (t_1,...,t_n)}{\part t_2},...,\frac{\part \Ga_\ga
(t_1,...,t_n)}{\part t_n}] dt_1...dt_n.\tag5 \endalign$$

{\smc Proof.}
Consider the $n$-parameter family of $(n-1)$-dimensional simplices
$\Delta(t_1,...,t_n)$  spanned by all
$n$-tuples of vertices $\Ga_\ga (t_1),\Ga_\ga (t_1,t_2),...,
\Ga_\ga (t_1,...,t_n),\; (t_1,...,t_n)\in
\Cal T^{n}$. By the above proposition a generic point on
  $\part \Cal {EH}_\ga$ is covered by this parametrization exactly
$n!$ times. Consider
now the element $d\omega$ of area
of $\part \Cal {EH}_\ga$ svept out by the small deformation of
$\Delta(t_1,...,t_n)$ when
$(t_1,...,t_n)$ are varying in $\Cal T^{n}$. Up to infinitesimals of higher
order
$d\omega$ coincides
with the area of the polytope  spanned by $\Ga_\ga (t_1),\Ga_\ga (t_1,t_2),...,
\Ga_\ga (t_1,...,t_n)$
and another $2^{n}$ vertices
of the form $\Ga_\ga (t_1+\epsilon_1
dt_1,t_2+\epsilon_2dt_2,...,t_n+
\epsilon_ndt_n)$ where each
$\epsilon_1,...,\epsilon_n$ is an arbitrary
$n$-tuple of zeros and ones. (We do not need to vary the vertices
$\Ga_\ga (t_1),\Ga_\ga (t_1,t_2)$ $ ,..., $
$\Ga_\ga (t_1,...,t_{n-1})$ since their small deformations belong to the
subspace $L_{t_1}$
of codimension 2 and
contribute to the higher order infinitesimals in $d\omega$.)  Now placing
the origin 0
 inside
$\Cal {EH}_\ga$ we can present the volume element $dV$ of the
 cone
over $d\omega$  as the volume of the polytope $\Pi$ spanned
by 0 and the above vertices. The last group of $2^{n}$ vertices of $\Pi$ form
the $n$-dimensional parallelepiped
up to the higher order infinitesimals. Let us split this  parallelepiped
into $n!$
$n$-dimensional
simplices and consider $n!$ $2n$-dimensional simplices spanned by each of
these simplices and the
rest of the vertices $0,\Ga_\ga (t_1),...,\Ga_\ga (t_1,...,t_{n})$. Volumes of
these $2n$-dimensional simplices coincide with each other up
to the higher order infinitesimals. Take one of these simplices, namely
the one spanned by
$0,\Ga_\ga (t_1),...,\Ga_\ga (t_1,...,t_{n-1}),
\Ga_\ga (t_1,...,t_n),$ $
\Ga_\ga (t_1+dt_1,t_2,...,t_n),$ $
\Ga_\ga (t_1,t_2+dt_2,...,t_n),...,$ $
\Ga_\ga (t_1,t_2,...,t_n+dt_n)$ with the same convention as above. Its
volume $Vol(Simp)$ with the appropriate orientation equals
$$\align Vol(Simp)=&\frac {1}{2n!}\det [\Ga _\ga(t_1),\Ga
_\ga(t_1,t_2),...,\Ga _\ga(t_1,...,t_n),\\
&\Ga _\ga(t_1+dt_1,...,t_n),  \Ga _\ga(t_1,...,t_n+dt_n)] . \tag6\endalign$$
Let us write $\Ga _\ga(t_1,...,t_i+dt_i,...,t_n)\approx
\Ga _\ga(t_1,...,t_i,...,t_n)+dt_i\frac{\part \Ga
_\ga(t_1,...,t_i,...,t_n)}{\part t_i}$ and substitute
 it in $(6)$.
The volume  of the whole polytope $\Pi$ is equal to $n!Vol(Simp)$.
Integrating against $\Cal T^{n}$ and taking into account
the fact that
 our parametrization covers $\part \Cal {EH}_\ga$ exactly
$n!$ times we get the
above integral formula
$(5)$.

\qed

\heading \S 3.  Appendix. Volume of $ \Cal {EH}_\ga$ for normalized
ellipse in $\bR^4$. \endheading

{\smc Proposition 3.1.} If $\ga$ is a normalized generalized ellipse in
$\bR^4$, i.e. $\ga=(\sin x, \cos x,$ $ 1/2\sin 2x,$ $ 1/2\cos x)$ then
$$Vol( \Cal {EH}_\ga)=\frac {9\pi^2\sqrt 3}{64}.$$

{\smc Proof.} Let $\gamma $ be a convex curve in $\Bbb R^4$. Then
$$
\Gamma (t_1,t_2) = \cases \gamma (t_1) + F(t_1,t_2)\dfrac{d\gamma }{dt}(t_1) +
   G(t_1,t_2)\dfrac{d^2\gamma }{dt^2}(t_1), &\text{if $t_2\ne t_1$}\\
                 \gamma (t_1), &\text{if $t_2=t_1$},\endcases
$$
where $F = \Delta_1/\Delta $, $G = \Delta _2/\Delta $, and
$$
\gather
\Delta _1 = \det (\gamma (t_2) - \gamma (t_1),\dfrac{d^2\gamma }{dt^2}(t_1),
\dfrac{d\gamma }{dt}(t_2),\dfrac{d^2\gamma }{dt^2}(t_2)),\\
\Delta _2 = \det (\dfrac{d\gamma }{dt}(t_1),\gamma (t_2) - \gamma (t_1),
\dfrac{d\gamma }{dt}(t_2),\dfrac{d^2\gamma }{dt^2}(t_2)),\
\Delta = \det (\dfrac{d\gamma }{dt}(t_1),\dfrac{d^2\gamma }{dt^2}(t_1),
\dfrac{d\gamma }{dt}(t_2),\dfrac{d^2\gamma }{dt^2}(t_2)).
\endgather
$$
Therefore
$$
Vol(\Cal {EH}_\gamma ) = \frac {1}{4!}\int _{\Cal T^2}G(F\dfrac{\partial G}
{\partial t_2} - G\dfrac{\partial F}{\partial t_2})W(t_1) dt_1dt_2
$$
where
$$
W(t) = \det (\gamma (t),\dfrac{d\gamma }{dt}(t),\dfrac{d^2\gamma }
{dt^2}(t),\dfrac{d^3\gamma }{dt^3}(t)).
$$
Integration by parts leads to
$$
Vol(\Cal {EH}_\gamma ) = - \frac {1}{16}\int _{\Cal T^2}
G^2\dfrac{\partial F}{\partial t_2}W(t_1) dt_1dt_2.
$$

Let $\gamma $ be the generalized Lissajoux curve
$$
\gamma (t) = (1/k\sin kt,1/k\cos kt,1/l\sin lt,1/l\cos lt)
$$
where $k,l$ are natural numbers, $k<l$. Then
$$
W(t)\equiv
\text{const} = \frac {1}{lk}(l^2-k^2)^2,
$$
$$
\align
\Delta _1& =
\frac {l^2-k^2}{lk}\left(k(1-\cos ls)\sin ks-
l(1-\cos ks)\sin ls\right)\\
&= \frac {l^2-k^2}{lk}\left(\cos (l-k)t - \cos (l+k)t \right)
\left((l+k)\sin (l-k)t - (l-k)\sin (l+k)t \right),\\
\Delta _2& =
\frac {1}{kl}\left(l^2(1+\cos ls)(1-\cos ks)+
k^2(1+\cos ks)(1-\cos ls)-2kl\sin ks\sin ls\right)\\
&= \frac {1}{kl}\left((l+k)\sin (l-k)t - (l-k)\sin (l+k)t \right)^2,\\
\endalign
$$
$$
\align
\Delta &=
2kl(1-\cos ks\cos ls)-(l^2+k^2)\sin ks\sin ls\\
&= (l+k)^2\sin ^2(l-k)t -(l-k)^2\sin ^2(l+k)t,
\endalign
$$
where $s=2t =t_2-t_1$. The mapping $\Gamma $ is
well-defined (i.e. $\Delta = 0$ iff $s = 0$ mod $2\pi$) if
and only if  $l-k=1$.
In the last case
$$
Vol(\Cal {EH}_\gamma ) = - \frac {\pi (2k+1)^2}{8k(k+1)}\int _{-\pi /2}^{\pi /2}
g^2\dfrac{df}{dt }dt
$$
where
$$
\gather
g(t )= \frac {(2k+1)\sin t - \sin (2k+1)t }
    {(k(k+1)\left((2k+1)\sin t + \sin (2k+1)t \right)},\\
f(t )=\frac {(2k+1)\left(\cos t - \cos (2k+1)t
      \right)}{k(k+1)\left((2k+1)\sin t + \sin (2k+1)t \right)},\\
\dfrac{df}{dt }= \frac {4(2k+1)\sin t \sin (2k+1)t }
    {\left((2k+1)\sin t + \sin (2k+1)t \right)^2}.
\endgather
$$

Note that $\sin (2k+1)t = h(t )\sin t $ where
$$
h(t )=\sum_{m=0}^k(-1)^m\binom {2k+1}{2m+1}
       (1-\cos ^2t )^m\cos ^{2(k-m)}t .
$$
Therefore
$$
Vol(\Cal {EH}_\gamma ) = - \frac {\pi (2k+1)^3}{2k^3(k+1)^3}\int _{-\pi
/2}^{\pi /2}
\frac{(2k+1-h)^2}{(2k+1+h)^4}hdt.
$$
Change of variable $z=\tan t $ leads to the formula
$$
Vol(\Cal {EH}_\gamma ) = - \frac {\pi (2k+1)^3}{2k^3(k+1)^3}\int _{-\infty }^
{+\infty }\frac{\left((2k+1)(1+z^2)^k-P_k(z^2)\right )^2}
{\left ((2k+1)(1+z^2)^k+P_k(z^2)\right )^4}P_k(z^2)(1+z^2)^{k-1}dz
$$
where
$$
P_k(x)=\sum_{m=0}^k(-1)^m\binom {2k+1}{2m+1}x^m
$$
is a polynomial of degree $k$.

Let us, for example, take $k=1$ (or $l=2$), i.e. consider the normalized
ellipse $\gamma =\Cal G_2$ in $\bR^4$. Then
$$
Vol(\Cal {EH}_\gamma )= \frac {27\pi }{16}\int _{-\infty }^
{+\infty }\frac{(z^2-3)z^4}{(3+z^2)^4}dz
= \frac {9\pi \sqrt 3}{16}(-2I_4+5I_3-4I_2+I_1)
$$
where
$$
I_m=\int _{-\infty }^{+\infty }\frac{dz}{(1+z^2)^m}
=\cases \pi ,&\text {if $m=1$}\\
\frac{1\cdot 3\cdot 5\cdot ...(2m-3)}
{2\cdot 4\cdot 6\cdot ...(2m-2)}\pi , &\text {if $m>1$}.
\endcases
$$
This gives
$$
Vol(\Cal {EH}_\gamma ) = \frac {9\pi ^2\sqrt 3}{64}.
$$

\qed

\heading \S 4. {Final remarks} \endheading

 Below we mention
some further questions.

1) Conjecture. For any 2 convex $\ga_1:S^1\to \rtn$ and
 $\ga_1:S^1\to \rtn$ and a Young diagram $\mu$ of area $n$ the pairs
$(\rtn, D_{\ga_1}(\mu))$ and $(\rtn, D_{\ga_2}(\mu))$ are homeomorphic. In
particular, all the standard discriminants $D_\ga$ are homeomorphic.

2) Find the projective analog of Schoenberg's inequality for the convex hull of a convex curve, i.e. consider
the Fubini-Study metric on $\bP^{2n}$ instead of the Eucledian structure.

3) Describe the essential set $\bold {Ess}_
\ga(\mu)$ for all $\mu$ and
obtain a
volume formula for all $\Cal {YH}_\ga(\mu)$.

4) Find volume formulas for $\Cal {YH}_\ga(\mu)$ for the case when the
target space is $P^n$ with the standard Fubini-Study metrics;

5) Find volume formulas for the dual domains $\Cal {YH}_\ga^*(\mu)$, i.e.
for the
convex hulls of different strata of $\bold {CSk}_\ga$.

6) Find volume formulas of other components of
the complement to $D_\ga(\mu)$, say, for the standard discriminant $D_\ga$ in
the Fubini-Stidy metric.

7) The final goal is to prove the conjecture formulated in the introduction
that fixing the
length of $\ga$ one gets the maximal volume of any $\Cal {YH}_\ga(\mu)$
only if $\ga$ is the generalized normalized ellipse.

\enddocument